\documentclass[11pt,reqno]{amsart}
\usepackage{amsmath,amssymb,geometry,color,tikz-cd}
\usepackage[backref,pagebackref,pdftex,hyperindex]{hyperref}
\usepackage[alphabetic]{amsrefs}
\usepackage{amssymb}% see geometry.pdf on how to lay out the page. There's lots.
\usepackage{bbm}
\usepackage{enumitem}
\usepackage{upgreek}

\newtheorem{proposition}{Proposition}[section]
\newtheorem{theorem}[proposition]{Theorem}
\newtheorem{lemma}[proposition]{Lemma}
\newtheorem{corollary}[proposition]{Corollary}
\newtheorem{conjecture}[proposition]{Conjecture}

\theoremstyle{definition}
\newtheorem{remark}[proposition]{Remark}
\newtheorem{definition}[proposition]{Definition}

\title{Abundance for large Kodaira dimension}
\author{Chuyu Zhou}
\address{\'Ecole Polytechnique F\'ed\'erale de Lausanne (EPFL), MA C3 615, Station 8, 1015 Lausanne, Switzerland}
\email{chuyu.zhou@epfl.ch}
\date{} % delete this line to display the current date

\newcommand{\Proj}{{\rm{Proj}}}

\newcommand{\Supp}{{\rm {Supp}}}

\newcommand{\bC}{\mathbb{C}}

\newcommand{\bN}{\mathbb{N}}

\newcommand{\bQ}{\mathbb{Q}}
\newcommand{\bR}{\mathbb{R}}

\newcommand{\bZ}{\mathbb{Z}}

\newcommand{\mO}{\mathcal{O}}

\newcommand{\tA}{\tilde{A}}

\begin{document}

\begin{abstract}
In this note, we apply the semi-ampleness criterion in Lemma \ref{lem: vertical} to prove many classical results in the study of abundance conjecture. As a corollary, we prove abundance for large Kodaira dimension depending only on \cite{BCHM10}.
\end{abstract}

\maketitle
\tableofcontents

\section{Introduction}
We work over $\bC$ throughout. The following conjecture lies in the central position in minimal model theory.

\begin{conjecture}\label{abundance conj}
Let $(X,\Delta)$ be a projective klt log pair of dimension $n\geq 2$ such that $\kappa(X,\Delta)\geq 0$, where $\kappa(X,\Delta)$ is the Kodaira dimension of $K_X+\Delta$. Suppose $K_X+\Delta$ is nef, then it is semi-ample.
\end{conjecture}

We first note here that nothing in this note is essentially new. However, as we found that the semi-ampleness criterion in Lemma \ref{lem: vertical} is powerful and convenient, the purpose of this note is to apply it to quickly prove many classical results in the study of Conjecture \ref{abundance conj} in a slightly different way, including Theorem \ref{large kod}, Theorem \ref{nef good}, and Theorem \ref{reduction}. We also derive an additive formula which precisely tells the difference between Iitaka dimension and numerical dimension, see theorem \ref{additive}. Finally, as a corollary, we prove the following widely known result depending only on \cite{BCHM10}. 

\begin{theorem}{\rm (=Corollary \ref{cor: large kod})}
Suppose Conjecture \ref{kod 0} is known up to dimension $l\in \bZ^+$. Let $(X,\Delta)$ be a minimal projective klt log pair of dimension $n\geq l+1$. If  $\kappa(X,\Delta)\geq n-l$, then $K_X+\Delta$ is semi-ample.
\end{theorem}

\noindent
\subsection*{Acknowledgement}
 The author is supported by grant European Research Council (ERC-804334).

\section{Preliminaries}

In this section, we list some preliminary results on very exceptional divisors. Readers may refer to \cite[Section 3]{Birkar12} for more details. We also refer to \cite{KM98, Laz04a} for the basic notion in birational geometry such as klt singularities, Iitaka dimension, and Kodaira dimension.

\begin{definition}{(\cite[Definition 3.2]{Sho03} or \cite[Definition 3.1]{Birkar12})}
Let $g: W\to Z$ be a contraction morphism  of projective normal varieties (i.e. $g_*\mO_W=\mO_Z$), $D$ an $\bR$-divisor on $W$, and $V\subset W$ a closed subset. We say that $V$ is verticle over $Z$ if $g(V)$ is a proper subset of $Z$. We say that $D$ is very exceptional over $Z$ if $D$ is verticle over $Z$ and for any prime divisor $P$ on $Z$ there is a prime divisor $Q$ on $W$ which is not a component of $D$ but $g(Q)=P$, i.e. over the generic point of $P$ we have $\Supp g^*P\nsubseteq \Supp D$.
\end{definition}

We fix some notation for this section. 

Let $X$ be a projective normal variety and $L$ a line bundle on $X$. Write $I(L)$ to be the Iitaka dimension of $L$. Suppose $I(L)\geq 1$ and $R(X,L)$ is finitely generated in degree one.  Let $f: W\to X$ be a resolution of $X$ and $f^*L\sim M+F$, where $M$ is globally generated and $F$ is the fixed part. Denote by $g: W\to Z$ the algebraic fiber space induced by $|kM|$ for a sufficiently divisible positive integer $k$. We note here that $\dim Z=I(L)$.

\begin{lemma}\label{lem: very exceptional}
Suppose $F$ is vertical over $Z$, then $F$ is very exceptional over $Z$.
\end{lemma}

\begin{proof}
This lemma is in fact implied by \cite[Lemma 3.2]{Birkar12}. We provide a proof here for readers' convenience. First note that there is an ample line bundle $\tA$ on $Z$ such that $kM\sim g^*\tA$. Since $R(X,L)$ is finitely generated by degree one, we know that $f^*rL\sim rM+rF$ where $rF$ is the fixed part for any positive integer $r$. Fix a positive integer $l$, we have the following
$$H^0(W, mkM+lF)\cong H^0(W, mg^*\tA+lF) \cong H^0(W, mg^*\tA)$$
for $m\gg1$. By projection formula, we see
$$H^0(Z, \mO_Z(m\tA)\otimes g_*\mO_W(lF)) \cong H^0(Z, \mO_Z(m\tA))$$
for $m\gg 1$. As there is an injection $\mO_Z\hookrightarrow g_*\mO_W(lF)$ and one can choose $m$ sufficiently large such that $m\tA$ is very ample, we have that $g_*\mO_W(lF)=\mO_Z=g_*\mO_W$.

Suppose $F$ is not very exceptional over $Z$, then there is a prime divisor $P$ on $Z$ such that if $Q$ is any prime divisor on $W$ such that $g(Q)=P$, then $Q$ is a component of $\Supp F$. Let $U$ be a smooth open subset of $Z$ such that $\mO_U\subsetneq \mO_U(P|_U)$, $P|_U$ is Cartier, and each component of $g^*P|_U$ maps onto $P|_U$. Let $V=g^{-1}U$. Then $\Supp g^*P\subset \Supp F|_V$, and $g^*P|_U\leq lF|_V$ for some $l>0$. However,
$$\mO_U\subsetneq \mO_U(P|_U)\subset g_*\mO_V(lF|_V)=\mO_U, $$
a contradiction.
\end{proof}

Let $S\to Z$ be a projective morphism of varieties and $M$ an $\bR$-Cartier divisor on $S$. We say that $M$ is nef on the very general curves of $S/Z$ if there is a countable union $\Lambda$ of proper closed subsets of $S$ such that $M.C\geq 0$ for any curve $C$ on $S$ contracted over $Z$ satisfying $C\nsubseteq \Lambda$.

\begin{lemma}\label{lem: negativity}{(\cite[Lemma 3.22]{Sho03}, \cite[Lemma 1.7]{Pro03}, or \cite[Lemma 3.3]{Birkar12})}
Let $g: W\to Z$ be a contraction of projective normal varieties. Let $D$ be an $\bR$-divisor on $W$ written as $D=D^+-D^-$ with $D^+, D^-\geq 0$ having no common components. Assume that $D^-$ is very exceptional over $Z$, and that for each component $S$ of $D^-$, $-D|_S$ is nef on the very general curves of $S/Z$. Then $D\geq 0$.
\end{lemma}

\begin{proof}
See the proof of \cite[Lemma 3.3]{Birkar12}.
\end{proof}

\begin{lemma}\label{lem: hyperplane}
Notation as in Lemma \ref{lem: very exceptional}. Suppose $F$ is very exceptional over $Z$ and $\dim g(F)>0$. Let $Z':=H$ be a very general hyperplane section of $Z$ and $W':=g^*H$. Denote by $g': W'\to Z'$ the restriction morphism, then $F|_{W'}$ is very exceptional over $Z'$.
\end{lemma}

\begin{proof}
The proof is contained in the proof of \cite[Lemma 3.3]{Birkar12}.
\end{proof}

\section{A criterion for semi-ampleness}

The following key lemma provides a criterion for semi-ampleness.

\begin{lemma}\label{lem: vertical}
Let $X$ be a projective normal variety and $L$ a nef line bundle on $X$. Suppose $I(L)\geq 1$ and $R(X,L)$ is finitely generated in degree one.   Let $f: W\to X$ be a resolution of $X$ and $f^*L\sim M+F$, where $M$ is globally generated and $F$ is the fixed part. Denote $g: W\to Z$ to be the algebraic fiber space induced by $|kM|$ for a sufficiently divisible positive integer $k$. Suppose  $F$ is vertical over $Z$, then $F=0$ and $L$ is semi-ample.
\end{lemma}

\begin{proof}

By Lemma \ref{lem: very exceptional}, $F$ is very exceptional. If $\dim g(F)=0$, then $F|_S\sim_\bQ (f^*L-g^*A)|_S$ is nef on the very general curves of $S$, where $S$ is any component of $F$. By Lemma \ref{lem: negativity}, we have $F\leq 0$. Thus $F=0$ and $f^*L\sim_\bQ g^*A$, which implies that $L$ is semi-ample.
If $\dim g(F)=l>0$, we choose $l$ very general hyperplane sections on $Z$, denoted by $H_1,..,H_{l}$,  and consider the induced morphism $$g': W':=g^*H_1\cap ...\cap g^*H_{l}\to Z':=H_1\cap...\cap H_{l} .$$
Suppose $F>0$, as $H_i, i=1,...,l$ are very general, we have $F|_{W'}>0$. By Lemma \ref{lem: hyperplane}, the divisor $F|_{W'}$ is very exceptional over $Z'$.  Note that we have the following $\bQ$-linear equivalence
$$f^*L|_{W'}\sim_\bQ g'^*(A|_{Z'})+F|_{W'}. $$
Since $F|_{W'}$ is very exceptional and $\dim g'(F|_{W'})=0$, 
$\{F|_{W'}\}|_S\sim_\bQ  \{f^*L|_{W'}- g'^*(A|_{Z'})\}|_S$
is nef on the very general curves of $S$, where $S$ is any component of $F|_{W'}$. By Lemma \ref{lem: negativity}, $F|_{W'}\leq 0$, which is a contradiction. Therefore $F=0$ and $f^*L\sim_\bQ g^*A$, which implies that $L$ is semi-ample. 
\end{proof}

As a corollary, we have the following well-known semi-ample result for large Iitaka dimension.

\begin{theorem}\label{large kod}
Let $X$ be a projective normal variety of dimension $n\geq 2$ and $L$ a nef line bundle on $X$. Suppose $I(L)\geq n-1$ and $R(X,L)$ is finitely generated, then $L$ is semi-ample.
\end{theorem}

\begin{proof}
Replace $L$ by a sufficiently large multiple, we may assume that $R(X,L)$ is finitely generated by degree one. Choose a resolution $f: W\to X$ such that $f^*L\sim M+F$, where $M$ is globally generated and $F$ is the fixed part. Let $g: W\to Z$ be the algebraic fiber space induced by $|kM|$ for a sufficiently large positive integer $k$, then we have a natural rational map $\phi: X\dashrightarrow Z$ where $Z=\Proj R(X,L)$, and there is an ample $\bQ$-line bundle $A$ on $Z$ such that 
$$f^*L\sim_\bQ g^*A+F .$$ 
By Lemma \ref{lem: vertical}, it suffices to show that $F$ is vertical over $Z$.
Suppose $I(L)=n$, then $\dim Z=n$ and $F$ is clearly vertical over $Z$. 
Suppose $I(L)=n-1$, then $\dim Z=n-1$.  If $F$ is not vertical over $Z$, then $F$ is clearly relatively big over $Z$. Thus there exist a relatively ample/$Z$ divisor $D$ and an effective $\bQ$-divisor $E$ on $W$ such that $F\sim_\bQ D+E+g^* B$, where $B$ is a $\bQ$-divisor on $Z$ which is not necessarily effective. Choose a sufficiently small rational number $0<\epsilon\ll 1$ such that 
\begin{enumerate}
\item $\frac{1}{2}g^*A+\epsilon D$ is ample on $W$;
\item $\frac{1}{2}A+\epsilon B$ is effective on $Z$.
\end{enumerate}
Therefore, we have the following expression:
$$f^*L=M+F\sim_\bQ \frac{1}{2}g^*A+\epsilon D+g^*(\frac{1}{2}A+\epsilon B)+\epsilon E+(1-\epsilon) F. $$
This expression indicates that $f^*L$ is big, which is a contradiction to the assumption that $I(L)=n-1$. The contradiction implies that $F$ is vertical over $Z$.
The proof is finished.
\end{proof}

\begin{corollary}
Let $X$ be a projective normal surface and $L$ a nef line bundle on $X$. Suppose $I(L)\geq 1$ and $R(X,L)$ is finitely generated, then $L$ is semi-ample.
\end{corollary}

\begin{definition}(\cite{Kawamata85})\label{def: num dim}
Let $X$ be a projective normal variety of dimension $n$ and $L$ a nef line bundle on $X$. The numerical dimension $\nu(L)$ is defined as follows:
$$\nu(L):=\max\{e\in \bN| \textit{$L^e.V\ne 0$ for some subvariety $V$ of dimension $e$}\}. $$
\end{definition}

As another application of Lemma \ref{lem: vertical}, we also prove the following well-known result for all positive Iitaka dimension (cf. \cite[Corollary 1]{MR97}).

\begin{theorem}\label{nef good}
Let $X$ be a projective normal variety of dimension $n\geq 2$ and $L$ a nef line bundle on $X$. Suppose $I(L)\geq 1$ and $R(X,L)$ is finitely generated. If $I(L)=\nu(L)$, then $L$ is semi-ample.
\end{theorem}

\begin{proof}
By Theorem \ref{large kod}, when $I(L)\geq n-1$, we even do not need the condition $I(L)=\nu(L)$. Thus it suffices to assume $1\leq I(L)\leq n-2$. We divide the proof into several steps.

\textbf{Step 0}. We use the same notation as in the proof of Theorem \ref{large kod}. Replace $L$ by a sufficiently large multiple, we may assume that $R(X,L)$ is finitely generated by degree one. Choose a resolution $f: W\to X$ such that $f^*L\sim M+F$, where $M$ is globally generated and $F$ is the fixed part. Let $g: W\to Z$ be the algebraic fiber space induced by $|kM|$ for a sufficiently large positive integer $k$, then we have a natural rational map $\phi: X\dashrightarrow Z$ where $Z=\Proj R(X,L)$, and there is an ample $\bQ$-line bundle $A$ on $Z$ such that 
$$f^*L\sim_\bQ g^*A+F .$$
By Lemma \ref{lem: vertical}, it suffices to show that $F$ is vertical over $Z$. Denote by $i:= n-I(L)-1\geq 1$ and suppose $F$ is not vertical over $Z$, we aim to derive a contradiction.

\textbf{Step 1}. Let $D_1$ be a relatively ample line bundle on $W$, and choose a sufficiently small rational number $0<\epsilon\ll 1$ such that $g^*A+\epsilon D_1$ is ample on $W$. Take a sufficiently divisible $m_1\in \bZ^+$ such that $|m_1(g^*A+\epsilon D_1)|$ is a very ample linear system on $W$. We choose a general member $H_1\in |m_1(g^*A+\epsilon D_1)|$, and it is clear that $F|_{H_1}$ still dominates $Z$. 
%Consider the following exact sequence:
%$$0\to \mO_W(f^*L-H_1)\to \mO_W(f^*L)\to \mO_{H_1}(f^*L)\to 0. $$
%By Kodaira vanishing, one can always assume $H^1(W, \mO_W(f^*L-H_1))=0$ by choosing $m_1$ sufficiently large. Thus the restriction map
%$$H^0(W, f^*L)\to H^0(H_1, f^*L|_{H_1}) $$ is surjective.
Note that we have the decomposition $f^*L|_{H_1}\sim M|_{H_1}+F|_{H_1}$. Replace $g: W\to Z$ with $g_1: H_1\to Z$ (which is the composition of $H_1\hookrightarrow W$ and $W\to Z$), and $f^*L$ (resp. $M, F$) with $f^*L|_{H_1}$ (resp. $M|_{H_1}, F|_{H_1}$).

\textbf{Step 2}. If $\dim F|_{H_1}> \dim Z$ (i.e. $i\geq 2$), we repeat the process as in step 1. Let $D_2$ be a relatively ample line bundle on $H_1$, and choose a sufficiently small rational number $0<\epsilon\ll 1$ such that $g_1^*A+\epsilon D_2$ is ample on $H_1$. Take a sufficiently divisible $m_2\in \bZ^+$ such that $|m_2(g_1^*A+\epsilon D_2)|$ is a very ample linear system on $H_1$. We choose a general member $H_2\in |m_2(g_1^*A+\epsilon D_2)|$, and it is clear that $F|_{H_2}$ still dominates $Z$. 
%Consider the following exact sequence:
%$$0\to \mO_{H_1}(f^*L|_{H_1}-H_2)\to \mO_{H_1}(f^*L|_{H_1})\to \mO_{H_2}(f^*L)\to 0. $$
%By Kodaira vanishing, one can always assume $H^1(H_1, \mO_{H_1}(f^*L|_{H_1}-H_2))=0$ by choosing $m_2$ sufficiently large. Thus the restriction map
%$$H^0(H_1, f^*L|_{H_1})\to H^0(H_2, f^*L|_{H_2}) $$ is surjective.
Note that we have the decomposition $f^*L|_{H_2}\sim M|_{H_2}+F|_{H_2}$. Replace $g_1: H_1\to Z$ with $g_2: H_2\to Z$ (which is the composition of $H_2\hookrightarrow H_1$ and $H_1\to Z$), and $f^*L|_{H_1}$ (resp. $M|_{H_1}, F|_{H_1}$) with $f^*L|_{H_2}$ (resp. $M|_{H_2}, F|_{H_2}$).

\textbf{Step 3}. If $\dim F|_{H_2}> \dim Z$, we repeat the above process. After $i$ steps, we have the morphism $g_i: H_i\to Z$ and the decomposition 
$$f^*L|_{H_i}\sim M|_{H_i}+F|_{H_i}\sim g_i^*A+F|_{H_i}.$$ 
Since $F|_{H_i}$ dominates $Z$ and $\dim F|_{H_i}=\dim Z$, we see that $F|_{H_i}$ is relatively big with respect to $g_i$. Thus there exist a relatively ample/$Z$ divisor $D$ and an effective $\bQ$-divisor $E$ on $H_i$ such that $F|_{H_i}\sim_\bQ D+E+g_i^* B$, where $B$ is a $\bQ$-divisor on $Z$ which is not necessarily effective. Choose a sufficiently small rational number $0<\epsilon'\ll 1$ such that 
\begin{enumerate}
\item $\frac{1}{2}g_i^*A+\epsilon' D$ is ample on $H_i$;
\item $\frac{1}{2}A+\epsilon' B$ is effective on $Z$.
\end{enumerate}
Therefore, we have the following expression:
$$f^*L|_{H_i}=M|_{H_i}+F|_{H_i}\sim_\bQ \frac{1}{2}g_i^*A+\epsilon' D+g_i^*(\frac{1}{2}A+\epsilon' B)+\epsilon' E+(1-\epsilon') F|_{H_i}. $$
This expression indicates that $f^*L|_{H_i}$ is big on $H_i$. As $\dim H_i=\dim Z+1=I(L)+1$, we see that $\nu(f^*L)\geq I(L)+1$, which implies that $\nu(L)\geq I(L)+1$. Contradiction.
\end{proof}

\section{Numerical dimension in family}

In this section, we fix $g: W\to Z$ to be a projective contraction morphism (i.e. $g_*\mO_W=\mO_Z$) of quasi-projective normal varieties with $\dim Z\geq 1$, and $F$ a $\bQ$-line bundle on $W$ such that $F_z$ is nef on $W_z$ for any closed point $z\in Z$, where $W_z$ is the fiber over $z\in Z$ and $F_z:=F|_{W_z}$.

\begin{lemma}\label{num dim}
Suppose $g: W\to Z$ is of relative dimension $d$, then there exists an open subset $U\subset Z$ such that the function $z\mapsto \nu(F_z)$ defined on $U$ is a constant function.
\end{lemma}

\begin{proof}
We choose $U$ to be the open subset of $Z$ such that $U$ is smooth and the morphism $W\times_Z U\to U$ is flat. Let $z_1, z_2\in U$ be two different closed points, we show that $\nu(F_{z_1})=\nu(F_{z_2})$. Suppose $\nu(F_{z_1})=e\in \bN$, then for any relatively ample line bundle $H$ on $W$, we have that 
$H^{d-e}.F^e. W_{z_1}\ne 0$. By \cite[Proposition 10.2]{Fulton84}, the following holds:
$$H^{d-e}.F^e. W_{z_1}=H^{d-e}.F^e. W_{z_2}\ne 0.$$
This implies that $\nu(F_{z_2})\geq \nu(F_{z_1})$. By symmetry, $\nu(F_{z_1})\geq \nu(F_{z_2})$, concluded.
\end{proof}

\section{Additive formula}

In this section, we fix $X$ to be a projective normal variety of dimension $n$ and $L$ a nef line bundle on $X$ with $I(L)\geq 1$.
We always assume that the graded ring $R(X,L)$ is finitely generated by degree one. Denote by $f: W\to X$ the resolution such that $f^*L\sim M+F$, where $M$ is globally generated and $F$ is the fixed part. Let $g: W\to Z$ be the algebraic fiber space induced by $|kM|$ for $k\gg 1$.
We will derive an additive formula which precisely reflects the difference between $I(L)$ and $\nu(L)$.
We note here that $F_z$ is nef on $W_z$ for general closed point $z\in Z$, since $M\sim_\bQ g^*A$ for some ample $\bQ$-line bundle $A$ on $Z$.  Combine the following additive formula and Theorem \ref{nef good}, one clearly sees that $L$ being semi-ample is equivalent to $\nu(F_z)=0$ for general $z\in Z$.

\begin{theorem}\label{additive}
Notation as above, we have the additive formula $\nu(L)=I(L)+\nu(F_z)$, where $z\in Z$ is a general closed point on $Z$.
\end{theorem}

\begin{proof}
By Theorem \ref{large kod}, it suffices to assume $1\leq I(L)\leq n-2$. By Lemma \ref{num dim}, $\nu(F_z)$ is constant for general closed point $z\in Z$. Denote by $i:=n-I(L)-1$ and $e:=\nu(F_z)$ for general closed point on $Z$.  Let $D$ be a relatively ample line bundle on $W$ and $A$ an ample $\bQ$-line bundle on $Z$ such that $M\sim_\bQ g^*A$, we choose a sufficiently small rational number $0<\epsilon\ll 1$ such that $g^*A+\epsilon D$ is ample on $W$. Consider the linear system $|m(g^*A+\epsilon D)|$ for a sufficiently large and divisible positive number $m$. Choose $i+1-e$ general elements in the linear system, denoted by $H_j, j=1,...,i+1-e$, then by \cite[Proposition 10.2]{Fulton84} we see that the intersection number $H_1...H_{i+1-e}.F^e.W_z$ is a non-zero constant for general closed point $z\in Z$. This means that $F_z|_{H_1\cap...\cap H_{i+1-e}\cap W_z}$ is big on $W_z|_{H_1\cap...\cap H_{i+1-e}}$ for general $z\in Z$, thus $F|_{H_1\cap...\cap H_{i+1-e}}$ is relatively big with respect to the morphism $H_1\cap...\cap H_{i+1-e}\to Z$ (which is the composition of $H_1\cap...\cap H_{i+1-e}\to W$ and $W\to Z$). Since 
$$f^*L|_{H_1\cap...\cap H_{i+1-e}}\sim f^*A|_{H_1\cap...\cap H_{i+1-e}}+F|_{H_1\cap...\cap H_{i+1-e}},$$
by the same explanation as in the proof of Theorem \ref{large kod}, the line bundle $f^*L|_{H_1\cap...\cap H_{i+1-e}}$ is big and nef on $H_1\cap...\cap H_{i+1-e}$. Note that 
$$\dim H_1\cap...\cap H_{i+1-e}=n-(i+1-e)=I(L)+e,$$ 
thus $\nu(L)=\nu(f^*L)\geq I(L)+\nu(F_z)$ for general $z\in Z$. The other direction follows from the next lemma.
\end{proof}

\begin{lemma}{\rm (\cite[Proposition V.2.7]{Nakayama04})}\label{easy dire}
Let $Y\to T$ be a projective contraction morphism from a smooth  projective variety such that $\dim T\geq 1$. Suppose $D$ is a pseudo-effective line bundle on $Y$, then 
$$\nu(D)\leq \dim T+\nu(D|_{Y_t}),$$  
where $Y_t$ is the fiber over a general closed point $t\in T$.
\end{lemma}

\begin{proof}
We only note here that the proof of \cite[Proposition V.2.7]{Nakayama04} uses another definition of numerical dimension which  defines for pseudo-effective line bundles, i.e. 
$$\nu(D)=\max\{e\in \bN|\textit{$\lim_{m\to \infty}\frac{\dim H^0(Y, mD+D')}{m^e}> 0$ for some ample $D'$ on Y}\}.$$
This definition coincides with Definition \ref{def: num dim} (cf. \cite[Proposition V.2.7(6)]{Nakayama04}) when $D$ is nef and the rest follows from the proof of \cite[Proposition V.2.7(9)]{Nakayama04}.
\end{proof}

\section{Reduction to  Kodaira dimension zero}

It is well-known that the following conjecture implies Conjecture \ref{abundance conj} (cf. \cite[Theorem 7.3]{Kawamata85}). In this section, we will prove this reduction by applying Lemma \ref{lem: vertical}.

\begin{conjecture}\label{kod 0}
Let $(X,\Delta)$ be a projective klt log pair of dimension $n\geq 2$ such that $\kappa(X,\Delta)=0$. Then we have $\nu(K_X+\Delta)=0$.
\end{conjecture}

\begin{remark}
We do not assume $(K_X+\Delta)$ to be minimal (i.e. $K_X+\Delta$ is nef) in the above conjecture. The numerical dimension $\nu(K_X+\Delta)$ of $K_X+\Delta$ is defined as in the proof of Lemma \ref{easy dire}.
\end{remark}

\begin{theorem}\label{reduction}
Conjecture \ref{kod 0} implies Conjecture \ref{abundance conj}.
\end{theorem}

\begin{proof}
Let $(X,\Delta)$ be a minimal projective klt pair of dimension $n\geq 2$. If $\kappa(K_X+\Delta)=0$, then $\nu(K_X+\Delta)=0$ by Conjecture \ref{kod 0}. This implies that $K_X+\Delta\equiv_\bQ 0$, thus $K_X+\Delta\sim_\bQ 0$ by \cite[Theorem 0.1]{Ambro05}. From now on, we assume $\kappa(X,\Delta)\geq 1$. By \cite{BCHM10}, we know that the graded ring $R(X,r(K_X+\Delta))$ is finitely generated, where $r$ is a positive integer such that $r(K_X+\Delta)$ is Cartier.
Denote by $L:=r(K_X+\Delta)$, we may assume that $r$ is sufficiently divisible such that $R(X,L)$ is finitely generated in degree one. Denote by $f: W\to X$ the log resolution of $(X,\Delta)$ such that $f^*L\sim M+F$, where $M$ is globally generated and $F$ is the fixed part. Let $g: W\to Z$ be the algebraic fiber space induced by $|kM|$ for $k\gg 1$. Write
$$K_W+f_*^{-1}\Delta+\Gamma=f^*(K_X+\Delta)+E, $$
where $(W, f_*^{-1}\Delta+\Gamma+E)$ is simple normal crossing, $\Gamma$ and $E$ are non negative divisors with no common components.
It is not hard to see that $g: W\to Z$ is also the Iitaka fibration with respect to the line bundle $f^*r(K_X+\Delta)+rE$, thus 
$$r(K_W+f_*^{-1}\Delta+\Gamma)|_{W_z}=(f^*r(K_X+\Delta)+rE)|_{W_z}\sim (F+rE)|_{W_z}$$
is of Iitaka dimension zero for general fiber $W_z$ (cf. \cite[Theorem 2.1.33]{Laz04a}). By adjunction, 
$$K_{W_z}+\Delta_{W_z}:= (K_W+f_*^{-1}\Delta+\Gamma)|_{W_z}$$
is a klt log pair of Kodaira dimension zero, thus by Conjecture \ref{kod 0}, 
$$\nu(K_{W_z}+\Delta_{W_z})=\nu((F+rE)|_{W_z})=0.$$
This also implies that $\nu(F|_{W_z})=0$ and thus $F$ is vertical over $Z$. By Lemma \ref{lem: vertical}, $L$ is semi-ample. The proof is finished.
\end{proof}

\begin{remark}
By the proof of Theorem \ref{reduction}, one sees that we can assume $(X,\Delta)$ to be simple normal crossing in Conjecture \ref{kod 0}, since $(W_z, \Delta_{W_z})$ is simple normal crossing for general $z\in Z$.
\end{remark}

\begin{corollary}\label{cor: large kod}
Suppose Conjecture \ref{kod 0} is known up to dimension $l\in \bZ^+$. Let $(X,\Delta)$ be a minimal projective klt log pair of dimension $n\geq l+1$. If  $\kappa(X,\Delta)\geq n-l$, then $K_X+\Delta$ is semi-ample.
\end{corollary}

\begin{proof}
We use the same notation as in the proof of Theorem \ref{reduction}. As $\kappa(X,\Delta)\geq n-l$, we see that $\dim W_z\leq l$. By our assumption and the proof of Theorem \ref{reduction}, $F$ is vertical and $K_X+\Delta$ is semi-ample.
\end{proof}

The following corollary is widely known and the proof does not depend on minimal model program once we are armed by \cite{BCHM10} and Conjecture \ref{kod 0} up to dimension 3.  One can also refer to \cite[Theorem 1.5]{Fili20} for a proof which depends on more MMP techniches.

\begin{corollary}
Let $(X,\Delta)$ be a minimal projective klt log pair of dimension $n\geq 4$. Suppose $\kappa(X,\Delta)\geq n-3$, then $K_X+\Delta$ is semi-ample.
\end{corollary}

\begin{proof}
The proof follows from the fact that Conjecture \ref{kod 0} is known up to dimension 3 (cf. \cite{KMM94}).
\end{proof}

\bibliography{reference.bib}
\end{document}